\documentclass[conference]{IEEEtran}
\usepackage{amsmath}
\usepackage{amssymb}
\usepackage{enumerate}

\usepackage{url}
\usepackage{algorithm}
\usepackage{algpseudocode}

\begin{document}
%
\title{Algorithms for structured matrix-vector product
of optimal bilinear complexity}

\author{\IEEEauthorblockN{Ke Ye}
\IEEEauthorblockA{Computational and Applied Mathematics Initiative\\
Department of Statistics\\
University of Chicago\\
Email: kye@galton.uchicago.edu}
\and
\IEEEauthorblockN{Lek-Heng Lim}
\IEEEauthorblockA{Computational and Applied Mathematics Initiative\\
Department of Statistics\\
University of Chicago\\
Email: lekheng@galton.uchicago.edu}
}


%


\maketitle

\begin{abstract}
We present explicit algorithms for computing structured matrix-vector products that are optimal in the sense of Strassen, i.e., using a provably minimum number of multiplications. These structures include Toeplitz/Hankel/circulant, symmetric, Toeplitz-plus-Hankel, sparse, and multilevel structures. The last category include \textsc{bttb}, \textsc{bhhb}, \textsc{bccb} but also any arbitrarily complicated nested structures built out of other structures.
\end{abstract}


%
\IEEEpeerreviewmaketitle

\section{Introduction}

Given a bilinear map $\beta: \mathbb{C}^m \times \mathbb{C}^n \to \mathbb{C}^p$, the \textit{bilinear complexity} \cite{Strass1,Strass2} of $\beta$ is the least number of multiplications needed to evaluate $\beta(x,y)$ for $x\in \mathbb{C}^m$ and $y\in \mathbb{C}^n$.  This notion of bilinear complexity is the standard measure of computational complexity for matrix inversion and matrix multiplication \cite{Strass0, LeGall,  Williams, Winograd}.

This article is an addendum to our work in \cite{YL} where we proposed a generalization of the \textit{Cohn--Umans method} \cite{CKSU,Cohn/Umans}  and used it to study the bilinear complexity of structured matrix-vector product. We did not derive any actual algorithms in \cite{YL}. The purpose of this present work is to provide explicit algorithms for structured matrix-vector product obtained by our generalized Cohn--Umans method in \cite{YL}. All algorithms in this paper have been shown to be the fastest possible in terms of bilinear complexity. The proofs may be found in \cite{YL} and involve determining the \textit{tensor ranks} of these structured matrix-vector products.

Here is a list of structured matrices discussed in this article:
\begin{itemize}
\item[\S\ref{circulant}] Circulant matrices.
\item[\S\ref{subsection:Toeplitz/Hankel}] Toeplitz/Hankel matrices.
\item[\S\ref{symmetric}] Symmetric matrices.
\item[\S\ref{Toeplitz-plus-Hankel}] Toeplitz-plus-Hankel matrices.
\item[\S\ref{sparse}] Sparse matrices.
\item[\S\ref{p-levels structured}] Multilevel structured matrices $A_1\otimes \cdots \otimes A_p$ where each $A_i$ is one of circulant, Toeplitz/Hankel, symmetric, Toeplitz-plus-Hankel, or sparse matrices.
\end{itemize}
The algorithms for circulant \cite{GL} and Toeplitz  \cite{BC} matrices are known but those for other structured matrices are new (as far as we know).  In particular, the multilevel structured matrices in \S\ref{p-levels structured}  include arbitrarily complicated nested structures, e.g., block \textsc{bccb} matrices whose blocks are Toeplitz-plus-Hankel, a $3$-level structure.

We analyze the bilinear complexities of all algorithms in \S\ref{sec:bc}.
Readers should bear in mind that bilinear complexity does not count scalar multiplications. For example, the bilinear map $\beta: \mathbb{C}^2 \times \mathbb{C}^2 \to \mathbb{C}$, $\beta((a,b),(c,d)) = (2a + b)(3c - d)$ has bilinear complexity one. For those familiar with tensor rank \cite{Strass}, the bilinear complexity of $\beta$ is just the tensor rank of the \textit{structure tensor} $\mu_\beta\in \mathbb{C}^2 \otimes \mathbb{C}^2\otimes \mathbb{C}^2$ corresponding to $\beta$ \cite{BCS,YL}.

\section{Circulant matrix}\label{circulant}
An $n\times n$ circulant matrix $A = (a_{ij})$ is a matrix with
\begin{align*}
a_{ij} &= a_{i+p,j+p},\quad  1\le i, j,i+p,j+p \le n,\\
a_{1j} &= a_{n+2-j,1}, \quad 2\le j \le n.
\end{align*}
The circulant matrix represented by $a = (a_1,\dots, a_n)\in \mathbb{C}^n$ is one whose first row is $a$. It is well-known \cite{GL} that the circulant matrix-vector product can be computed by Fourier transform. We restate this algorithm for completeness. Let $\omega_k = e^{2k\pi i/n}$, $k=0,\dots, n-1$ and define the Fourier matrix 
\begin{equation}\label{eq:F}
W = \begin{bmatrix}
1 & 1 & \cdots & 1\\
1 & \omega_1 & \cdots & \omega_{n-1}\\
\vdots & \vdots & \ddots & \vdots\\
1 & \omega_1^{n-1} & \cdots & \omega_{n-1}^{n-1}
\end{bmatrix}.
\end{equation}
\begin{algorithm}
  \caption{Circulant matrix-vector product}
    \label{algorithm:circulant}
  \begin{algorithmic}[1]
    \State Represent the circulant matrix $A$ by $a = (a_1, a_2, \dots, a_n)^{\mathsf{T}}$ and the column vector by $v = (v_1,v_2,\dots,v_n)^{\mathsf{T}}$.
    \State Compute $Wa$ and represent it by $(\tilde{a}_1,\dots, \tilde{a}_n)^{\mathsf{T}}$.
    \State Compute $nW^{-1}v$ and represent it by $(\tilde{v}_1,\dots, \tilde{v}_n)^{\mathsf{T}}$.
    \State Compute $\tilde{z} = (\tilde{a}_1\tilde{v}_1, \dots, \tilde{a}_n\tilde{v}_n)^{\mathsf{T}}$.
    \State Compute $z = W \tilde{z}$, which is the product of $A$ and $v$.
  \end{algorithmic}
\end{algorithm}

\section{Toeplitz/Hankel matrix}\label{subsection:Toeplitz/Hankel}

An $n\times n$ Toeplitz matrix $A = (a_{ij})$ is a matrix with
\[
a_{ij} = a_{i+p,j+p}, \quad 1\le i,j, i+p,j+p \le n.
\]
We represent an $n\times n$ Toeplitz matrix $A = (a_{ij})$ by $(a_1, a_2, \dots, a_{2n-1})\in \mathbb{C}^{2n-1}$  
\[
a_{ij} = a_{j - i + n }. 
\]
Every $n\times n$ Toeplitz matrix $A$ may be regarded as a block of some $2n\times 2n$ circulant matrix $C$ whose first row is $(a_n,\dots, a_{2n-1},b,a_1,\dots, a_{n-1})$ and $b\in \mathbb{C}$ is arbitrary. Using this embedding, we obtain Algorithm~\ref{algorithm:Toeplitz} for Toeplitz matrix-vector product \cite{BC, YL} .
\begin{algorithm}
  \caption{Toeplitz matrix-vector product}
    \label{algorithm:Toeplitz}
  \begin{algorithmic}[1]
\State Express the Toeplitz matrix $A$ as $(a_1,\dots, a_{2n-1})$ and the vector as $v = (v_1,\dots, v_n)^{\mathsf{T}}$.
\State Compute $b = -\sum_{i=1}^{2n-1} a_i$.
\State Construct $c = (a_n,\dots, a_{2n-1},b,a_1,\dots, a_{n-1}) \in \mathbb{C}^{2n}$.
\State Construct $\tilde{v} = (v_1,\dots, v_n, 0, \dots, 0)^{\mathsf{T}}\in \mathbb{C}^{2n}$.
\State Compute the product $\tilde{z} = (z_1,\dots, z_{2n})^{\mathsf{T}}$ of the circulant matrix determined by $c$ with $\tilde{v}$ by Algorithm~\ref{algorithm:circulant}. 
\State $z = (z_1,\dots, z_n)^{\mathsf{T}}$ is the product of $A$ and $v$.
  \end{algorithmic}
\end{algorithm}

An $n\times n$ matrix $H$ is called a Hankel matrix if $JH$ is a Toeplitz matrix where 
\begin{equation}\label{eq:J}
J =\begin{bmatrix}
0 & 0 &\cdots & 0 & 1\\
0 & 0& \cdots & 1 & 0\\
\vdots & \vdots & \ddots & \vdots & \vdots\\
0 & 1& \cdots & 0 & 0\\
1 & 0& \cdots & 0 & 0\\
\end{bmatrix}.
\end{equation}
We represent an $n \times n$ Hankel matrix $H = (h_{ij})$ as $(h_1,h_2,\dots, h_{2n-1}) \in \mathbb{C}^{2n - 1}$ where 
\[
h_{ij} = h_{2n+1-i-j}, \quad 1\le i, j\le n.
\]
Algorithm~\ref{algorithm:Hankel} computes the product of a Hankel matrix and a column vector $v$.
\begin{algorithm}
  \caption{Hankel matrix-vector product}
    \label{algorithm:Hankel}
  \begin{algorithmic}[1]
  \State  Express $T = (h_1,h_2,\dots, h_{2n-1})$.
    \State Apply Algorithm~\ref{algorithm:Toeplitz} to the Toeplitz matrix represented by $T$ and $v$ to obtain $(z_1,\dots, z_n)$.
    \State $(z_n, z_{n-1},\dots, z_1)$ is the product of $H$ and $v$.
  \end{algorithmic}
\end{algorithm}

\section{Symmetric matrix}\label{symmetric}
Algorithm~\ref{algorithm:Symmetric} computes the product of a symmetric matrix $S=(s_{ij})$ where $s_{ij} = s_{ji}$ and a column vector $v$. We represent a symmetric matrix $s = (s_{ij})$ as $(s_1,\dots, s_{N}) \in \mathbb{C}^N$ where $N = \binom{n+1}{2}$ and the index of $s_k$ is
\[
k = (i-1)n - \binom{i-1}{2} + j, \quad 1\le i \le j \le n.
\]
\begin{algorithm}
  \caption{Symmetric matrix-vector product}
    \label{algorithm:Symmetric}
  \begin{algorithmic}[1]
  \State $S$ is an $n\times n$ symmetric matrix. Set $S_1 = S$. Set $m = \lceil n/2 \rceil$.  Set $v_1 = v$ and $z=0$.
        \For {$k =1,\dots, m$}
    \State Construct Hankel matrix $H_k$ determined by first row and last column of $S_k$.
    \State Compute $w_ k = H_k v_k$ by Algorithm~\ref{algorithm:Hankel}.
    \State Update $z = z + w_k$.
    \State Construct $S_{k+1}$ by deleting first and last columns and  first and last rows of $S_k- H_k$.
     \State Construct $v_{k+1}$ by deleting first and last entry of $v_k$.
       \EndFor
      \State $z = (z_1,\dots, z_n)^{\mathsf{T}}$ is the product of $S$ and $v$.
  \end{algorithmic}
\end{algorithm}

\section{Toeplitz-plus-Hankel matrix}\label{Toeplitz-plus-Hankel}
An $n\times n$ Toeplitz-plus-Hankel matrix is a matrix which can be expressed as the sum of an $n\times n$ Hankel matrix and an $n\times n$ Toeplitz matrix. If $X$ is an $n\times n$ Toeplitz-plus-Hankel matrix and 
\[
X = H + T
\] 
for some Hankel matrix $H$ and some Toeplitz matrix $T$, then for any $a\in \mathbb{C}$ we have a decomposition of $X$ into the sum of a Hankel matrix $H + aE$ and a Toeplitz matrix $T - aE$ where $E$ is the $n\times n$ matrix with all entries equal to one.

\begin{algorithm}
  \caption{Toeplitz-plus-Hankel matrix-vector product}
    \label{algorithm:Toeplitz-plus-Hankel}
  \begin{algorithmic}[1]
  \State Express $X$ as $H + T$ with Hankel matrix $H$ and Toeplitz matrix $T$. 
  \State Express $T$ as $(t_1,\dots, t_{2n-1})$ and $H$ as $(h_1,\dots,h_{2n-1})$.
  \State Compute $b = -\sum_{j=1}^{2n-1} t_j$.
  \State Compute $a\in \mathbb{C}$ as
  \[
a =\frac{\sum_{j=0}^{n-1} \omega^{j} t_{n+j} + \omega^n b + \sum_{j=1}^{n-1} \omega^{n+j} t_j}{2n}
   \] 
where $\omega= e^{k \pi i/n}$.
  \State Update $H = H + aE$ and $T = T - aE$.
  \State Compute $z_H = H v$ by Algorithm~\ref{algorithm:Hankel} and $z_T =Tv$ by Algorithm~\ref{algorithm:Toeplitz}, respectively.
  \State Compute $z = z_H + z_T$, which is the product of $X$ and $v$.
  \end{algorithmic}
\end{algorithm}

\section{Sparse matrix}\label{sparse}
An $n\times n$ sparse matrix $A = (a_{ij})$ with sparsity pattern $\Omega \subseteq \{1,\dots, n\} \times \{1,\dots, n\}$ is one where 
\[
a_{ij} = 0 \quad \text{for all}\; (i,j)\in \Omega.
\]
For example, an upper triangular matrix is a sparse matrix with sparsity pattern $\Omega=\{(i,j):1\le i\le j \le n\}$. For sparse matrices associated with $\Omega$, the matrix-vector product has optimal bilinear complexity $\#\Omega$ realized by the usual matrix-vector product algorithm \cite{YL}.

\section{Multilevel structured matrix}\label{p-levels structured}
Let $A=(a_{ij})\in \mathbb{C}^{n\times n}$ and $B=(b_{ij}) \in \mathbb{C}^{m\times m}$. The Kronecker product \cite{VL} of $A$ and $B$ is defined as
\[
A\circledast B = (a_{ij} B) \in \mathbb{C}^{mn \times mn},
\] 
i.e., $A\circledast B$ is an $m\times m$ block matrix whose $(i,j)$th block is the $n\times n$ matrix $a_{ij}B$. We may iterate the definition to obtain a $p$ levels matrix $A =A_1\circledast \cdots \circledast A_p$. In particular, if $A_1,\dots, A_p$ are structured matrices (circulant, Toeplitz, Hankel, symmetric and Toeplitz-plus-Hankel), then $A$ is called a $p$ levels structured matrix. 

Let $X_1\subseteq \mathbb{C}^{n_1 \times n_1},\dots,X_p\subseteq \mathbb{C}^{n_p \times n_p}$ be subspaces of structured matrices. Then $X_1\circledast \cdots \circledast X_p \subseteq \mathbb{C}^{n_1\cdots n_p \times n_1 \cdots n_p}$ is the set of all $p$ levels structured matrices $A_1\circledast \cdots \circledast A_p$ where $A_1 \in X_1,\dots,A_p \in X_p$.

Algorithms~\ref{algorithm:p-levels structured 1}--\ref{algorithm:p-levels structured 6} are based on the following idea. Let $\beta_i: X_i \times \mathbb{C}^{n_i} \to \mathbb{C}^{n_i}$ be the bilinear map defined by the matrix-vector product for matrices in $X_i$. Assume that the bilinear complexity of $\beta_i$ is $r_i$. Then the structural tensor \cite{YL} $\mu_{\beta_i} \in X_i^* \otimes (\mathbb{C}^{n_i})^* \otimes \mathbb{C}^{n_i}$ of $\beta_i$ has a tensor decomposition
\[
\mu_{\beta_i} = \sum_{j=1}^{r_i} u_j\otimes v_j\otimes w_j.
\]
The bilinear map $\beta: (X_1\circledast \cdots \circledast X_p)\times \mathbb{C}^{n_1\cdots n_p} \to \mathbb{C}^{n_1\cdots n_p}$, defined by the $p$ levels structured matrix-vector product, has structural tensor $\mu_\beta = \mu_{\beta_1}\otimes \cdots \otimes \mu_{\beta_p}$. In \cite{YL} we showed that if $X_i$ is Toeplitz, Hankel, symmetric, or Toeplitz-plus-Hankel, the bilinear complexity is equal to the dimension of $X_i$ and we obtain a machinery to decompose $\mu_{\beta_i}$ explicitly. Essentially, Algorithms~\ref{algorithm:p-levels structured 1}--\ref{algorithm:p-levels structured 6} are obtained from the tensor decompositions of structural tensors.

\subsection{Illustrative example}\label{example}
As an example, let us consider the case where $p=2$ and $A,B$ are $2\times 2$ circulant matrices. This gives a block-circulant-circulant-block or \textsc{bccb} matrix.  We set 
\[
A = 
\begin{bmatrix}
a & b\\
b & a
\end{bmatrix}, \quad
B =\begin{bmatrix}
c & d\\
d & c
\end{bmatrix},
\]
and 
\[
v = (x,y,z,w)^{\mathsf{T}} = \begin{bmatrix}
x \\ y
\end{bmatrix}\circledast \begin{bmatrix}
1 \\ 0
\end{bmatrix} +  \begin{bmatrix}
z \\ w
\end{bmatrix}\circledast \begin{bmatrix}
0 \\ 1
\end{bmatrix}.
\]
We want to compute the product of $A\circledast B$ and $v$. By definition we have
\[
A\circledast B = 
\begin{bmatrix}
aB & bB \\
bB & aB\\
\end{bmatrix}
 = 
\begin{bmatrix}
ac & ad & bc & bd \\
ad & ac & bd & bc\\
bc & bd & ac & ad\\
bd & bc & ad & ac
\end{bmatrix}
, 
\]
and 
\begin{align*}
(A\circledast B) v 
&= 
\begin{bmatrix}
a(\xi_1 + \xi_2) + b(\eta_1 + \eta_2)\\
a(\xi_1 - \xi_2) + b(\eta_1 - \eta_2)\\
b(\xi_1 + \xi_2) + a(\eta_1 + \eta_2)\\
b(\xi_1 - \xi_2) + a(\eta_1 - \eta_2)\\
\end{bmatrix}
,
\end{align*}
where 
\begin{align*}
\xi_1 &= \frac{1}{2}((cx + dy) + (dx + cy)), \\
\xi_2 &= \frac{1}{2}((cx+dy)-(dx + cy)),\\
\eta_1 &= \frac{1}{2}((cz + dw) + (dz + cw)),\\
\eta_2 &= \frac{1}{2}((cz+dw)-(dz + cw)).
\end{align*}
Observe that 
\begin{align*}
\begin{bmatrix}
a(\xi_1 + \xi_2) + b(\eta_1 + \eta_2)\\ 
b(\xi_1 + \xi_2) + a(\eta_1 + \eta_2)
\end{bmatrix}
&= \frac{1}{2} 
\begin{bmatrix}
\alpha + \beta \\
\alpha - \beta
\end{bmatrix}
,
\end{align*}
where
\begin{align*}
\alpha &= (a+b)[(\xi_1 + \xi_2) + (\eta_1 + \eta_2)] \\
&= (a+b) [(\xi_1+\eta_1) + (\xi_2+ \eta_2)], \\
\beta &= (a-b)[(\xi_1 + \xi_2) - (\eta_1 + \eta_2)]\\
&=(a-b) [(\xi_1-\eta_1) + (\xi_2 - \eta_2)].
\end{align*}
Similarly, we have
\begin{align*}
\begin{bmatrix}
a(\xi_1 - \xi_2) + b(\eta_1 - \eta_2)\\ 
b(\xi_1 - \xi_2) + a(\eta_1 - \eta_2)
\end{bmatrix}
&= \frac{1}{2} 
\begin{bmatrix}
\gamma + \tau \\
\gamma - \tau
\end{bmatrix}
,
\end{align*}
where
\begin{align*}
\gamma &= (a+b)[(\xi_1 - \xi_2) + (\eta_1 - \eta_2)] \\
& = (a+b)[(\xi_1 + \eta_1) - (\xi_2 + \eta_2)], \\
\tau &= (a-b)[(\xi_1 - \xi_2) - (\eta_1 - \eta_2)] \\
&= (a-b)[(\xi_1 - \eta_1) - (\xi_2 - \eta_2)].
\end{align*} 
Lastly, we observe that 
\begin{align*}
\xi_1 + \eta_1 &= \frac{1}{2} (c+d)[(x + y) + (z + w)],\\
\xi_1 - \eta_1 &= \frac{1}{2} (c+d)[(x + y) - (z + w)],\\
\xi_2 + \eta_2 &= \frac{1}{2} (c-d)[(x - y) + (z - w)],\\
\xi_2 - \eta_2 &= \frac{1}{2} (c-d)[(x - y) - (z - w)].\\
\end{align*}
By above computations, we see that one may compute $(A\circledast B) v$ using four multiplications, i.e., it is sufficient to compute 
\begin{align*}
w_{11} &=(a+b) (c+d)[(x + y) + (z + w)], \\
w_{12} &=(a+b) (c-d)[(x - y) + (z - w)], \\
w_{21} &=(a-b) (c+d)[(x + y) - (z + w)], \\
w_{22} &=(a-b) (c-d)[(x - y) - (z - w)].
\end{align*}
Note that since the entries of $A\circledast B$ are given as inputs, evaluating terms like $(a+b)(c+d) = ac + ad +bc +bd$ does not cost any multiplication (as we already have $ac,ad,bc,bd$ as inputs).

\subsection{General case}
We now generalize the above calculations to obtain an algorithm for $p$ levels structured matrix-vector product. In order to treat all cases at one go, our presentation in this section is slightly more abstract. Given a $p$ levels structured matrix $B \in X_1 \circledast \dots \circledast X_p$ and a vector $v$ of appropriate size, our algorithm, when applied to $B$ and $v$, takes the form:
\[
(B,v) \xrightarrow{(\varphi,\psi)} (b',v') \xrightarrow{m} m(b',v') \xrightarrow{\vartheta} Bv,
\]
where $\varphi$ is a linear map sending $B$ to a vector $b'$, $\psi$ is a linear map sending $v$ to a vector $v'$, $m$ is pointwise multiplication, and $\vartheta$ is another linear map sending $m(b',v')$ to $Bv$.  $\varphi$, $\psi$, and $\vartheta$ depend only on the structure of $B$ (i.e., on $X_1,\dots,X_p$) but not on the values of $B$ and $v$. For any given structure, we can represent the linear maps $\varphi$, $\psi$, and $\vartheta$ concretely as matrices.

We will present the  algorithms for $p$ levels structured matrix-vector product \textit{inductively}, by calling the corresponding $p-1$ levels algorithms. Also, they will be built upon Algorithms~\ref{algorithm:Toeplitz}, \ref{algorithm:Hankel}, \ref{algorithm:Symmetric}, and \ref{algorithm:Toeplitz-plus-Hankel} for the relevant structured matrix-vector product.

Suppose we have algorithms for $p-1$ levels structured matrix-vector product, i.e., we may evaluate the linear maps $\varphi$, $\psi$, and $\vartheta$ for any $p-1$ levels structured matrix.
Given a $p$ levels structured matrix $A_1\circledast \cdots \circledast A_p$ and a column vector $v$ of size $N =\prod_{i=1}^p n_i$, we write $A_1 \circledast \cdots \circledast A_p$ as $A \circledast B$ where $A = A_1$ and $B = A_2\circledast \cdots \circledast A_p$. Set $N_1$ to be $N/n_1$.

Let $A$ be a circulant matrix. Let $\omega_k= e^{2k\pi i /n}$, $k=0,1,\dots, n-1$ be the $n$th roots of unity and let $W = (\omega_k^j)_{j,k=0}^{n-1}$ be the Fourier matrix  in \eqref{eq:F}. We have Algorithm~\ref{algorithm:p-levels structured 1}.
\begin{algorithm}
  \caption{$p$ levels circulant matrix-vector product}
    \label{algorithm:p-levels structured 1}
  \begin{algorithmic}[1]
\State Express $A$ by a column vector $a = (a_1,\dots, a_{n_1})^{\mathsf{T}}$ and express $v$ by a column vector 
\begin{multline*}
v=(v_{1,1},\dots, v_{1,N_1}, v_{2,1},\dots, v_{2,N_1},\dots,\\
 v_{n_1,1},\dots, v_{n_1,N_1})^{\mathsf{T}}.
\end{multline*}
\State Express $\varphi$ as $(\varphi_1,\dots, \varphi_r)^{\mathsf{T}}$ where $\varphi_j$ is a linear functional on $X_2\circledast \cdots \circledast X_p$ and $r = \prod_{i=2}^p \dim(X_i)$. 
\State Express $\psi$ as $(\psi_1,\dots, \psi_{r})$ where $\psi_j$ is a linear functional on $\mathbb{C}^{N_1}$.
\State Compute $\tilde{a} = W a$ and denote it by $(\tilde{a}_{1},\dots,\tilde{a}_{n_1})^{\mathsf{T}}$.
\State Denote $v_ i = (v_{i,1},\dots, v_{i,N_1})^{\mathsf{T}},i=1,\dots, n_1$.
\For{$s=1,\dots,n_1$}
\For{$t=1,\dots, r$}
\State Compute 
\[
w_{st} =  \tilde{a}_s\varphi_t(B)\sum_{k=1}^{n_1} \omega_{k-1}^{s-1} \psi_t(v_k).
\]
\EndFor
\EndFor
\State Represent $(w_{st})$ as a column vector 
\begin{multline*}
w= (w_{11},\dots, w_{1,r},w_{21},\dots,w_{2,r},\dots,\\
w_{n_1,1},\dots,w_{n_1,r})^{\mathsf{T}}.
\end{multline*}
\State Compute $ (W \circledast \vartheta) w$, which is the product  $(A\circledast B)v$. 
\end{algorithmic}
\end{algorithm}

If we apply Algorithm~\ref{algorithm:p-levels structured 1} to the case where $A,B$ are $2\times 2$ circulant matrices, we obtain $w_{11},w_{12},w_{21},w_{22}$ as in Section~\ref{example}. To  compute the product of $A \circledast B$ and $v$, we express $A$ as $(a,b)^{\mathsf{T}}$, $B$ as $(c,d)^{\mathsf{T}}$, and $v$ as $(x,y,z,w)^{\mathsf{T}}$. Hence $v_1 = (x,y)^{\mathsf{T}}$ and $v_2 = (z,w)^{\mathsf{T}}$. By Algorithm~\ref{algorithm:circulant} the linear map $\varphi = (\varphi_1,\varphi_2)^{\mathsf{T}}$ is given by $\varphi_1((\alpha,\beta)^{\mathsf{T}}) = \alpha + \beta$ and $\varphi_2((\alpha,\beta)^{\mathsf{T}}) = \alpha - \beta$, and $\psi$ is the map given by $\psi_1((\alpha,\beta)^{\mathsf{T}}) = \alpha+\beta$ and $\psi_2((\alpha,\beta)^{\mathsf{T}}) = \alpha-\beta$, where $(\alpha,\beta)^{\mathsf{T}}$ is any column vector of size two. Lastly, the linear map $\vartheta$ is given by left multiplication by $\left[\begin{smallmatrix}
1 & 1\\
1 & -1
\end{smallmatrix}\right]$.

Let $A$ be a Toeplitz matrix. As before, there exists a circulant matrix $C$ of the form
\[
C =
\begin{bmatrix}
A & A' \\
A' & A
\end{bmatrix},
\]
and 
\[
(C \circledast B) \begin{bmatrix}
v\\
0
\end{bmatrix} = \begin{bmatrix}
(A\circledast B) v \\
(A' \circledast B) v
\end{bmatrix}.
\]
Hence to compute $(A\circledast B)v$, it suffices to compute $(C\circledast B)\left[\begin{smallmatrix}
v\\
0
\end{smallmatrix}\right]$ and this can be done using Algorithm~\ref{algorithm:p-levels structured 1}. We obtain Algorithm~\ref{algorithm:p-levels structured 2}.

\begin{algorithm}
  \caption{$p$ levels Toeplitz matrix-vector product}
    \label{algorithm:p-levels structured 2}
  \begin{algorithmic}[1]
  \State Express $A$ as a vector $a = (a_1,\dots, a_{2n-1})$ and $v$ as $(v_1,\dots, v_N)^{\mathsf{T}}$.
  \State Compute $b = -\sum_{i=1}^{2n_1-1} a_i$.
  \State Construct  $c = (a_n,\dots, a_{2n_1-1},b,a_{1},\dots, a_{n_1-1}) \in \mathbb{C}^{2n_1}$ representing a $2n_1\times 2n_1$ circulant matrix $C$.
  \State Construct $\tilde{v} = (v_1,\dots, v_N, 0 ,\dots, 0)\in \mathbb{C}^{2N}$.
  \State Compute $\tilde{z} = (C\circledast B)\tilde{v}$ by Algorithm~\ref{algorithm:p-levels structured 1} and express $\tilde{z}$ as $(z_1,\dots, z_{2N})^{\mathsf{T}}$.
  \State $(z_1,\dots, z_N)^\mathsf{T}$ is the product of $(A\circledast B)$ and $v$.
\end{algorithmic}
\end{algorithm}

Now for square Hankel matrices $A$ and $B$ we observe that
\[
JA\otimes B =( J\otimes I) (A\otimes B),
\]
where $J$ is the matrix  in \eqref{eq:J}. Algorithm~\ref{algorithm:p-levels structured 3} follows.
\begin{algorithm}
  \caption{$p$ levels Hankel matrix-vector product}
    \label{algorithm:p-levels structured 3}
  \begin{algorithmic}[1]
 \State Compute the $Z = (JA\otimes B) v$ by Algorithm~\ref{algorithm:p-levels structured 2}.
 \State Compute $z = (J\otimes I)Z$ and $z$ is $(A\otimes B)v$.
\end{algorithmic}
\end{algorithm}

The algorithms for $p$ levels symmetric matrix (Algorithm~\ref{algorithm:p-levels structured 4}),  $p$ levels Toeplitz-plus-Hankel matrix  (Algorithm~\ref{algorithm:p-levels structured 5}),  $p$ levels sparse matrix  (Algorithm~\ref{algorithm:p-levels structured 6}) are obtained via similar considerations.
\begin{algorithm}
  \caption{$p$ levels symmetric matrix-vector product}
    \label{algorithm:p-levels structured 4}
  \begin{algorithmic}[1]
  \State $A$ is an $n_1\times n_1$ symmetric matrix. Compute $m = \lceil n_1/2\rceil$.  Set $v_1 = v$ and $z=0\in\mathbb{C}^{N}$.
        \For {$k =1,\dots, m$}
    \State Construct $H_k$ determined by first row and last column of $A_k$.
    \State Compute $w_ k = (H_k\circledast B) v_k$ by Algorithm~\ref{algorithm:p-levels structured 3}.
    \State Update $z = z + w_k$.
    \State Construct $A_{k+1}$ by deleting first and last columns and first and last rows of $A_k- H_k$.
     \State Construct $v_{k+1}$ by deleting first $N_1$ and last $N_1$ entries of $v_k$.
       \EndFor
      \State $z = (z_1,\dots, z_N)^{\mathsf{T}}$ is the product of $S$ and $v$.
  \end{algorithmic}
\end{algorithm}

\begin{algorithm}
  \caption{$p$ levels Toeplitz-plus-Hankel matrix-vector product}
    \label{algorithm:p-levels structured 5}
  \begin{algorithmic}[1]
  \State Express $A$ as $H + T$ with Hankel matrix $H$ and Toeplitz matrix $T$.
  \State Express $T$ as $(t_1,\dots, t_{2n_1-1})$ and $H$ as $(h_1,\dots,h_{2n_1-1})$.
  \State Compute $b = -\sum_{j=1}^{2n_1-1} t_j$.
  \State Find $a\in \mathbb{C}$ such that 
  \[
a =\frac{\sum_{j=0}^{n_1-1} \omega_1^{j} t_{n_1+j} + \omega_1^{n_1} b + \sum_{j=1}^{n_1-1} \omega_1^{n_1+j} t_j}{2n_1}
   \] 
where $\omega_1= e^{k \pi i/n_1}$.
  \State Update $H = H + aE$ and $T = T - aE$.
  \State Compute $z_H = (H\circledast B) v$ by Algorithm~\ref{algorithm:p-levels structured 3} and $z_T =(T\circledast B)v$ by Algorithm~\ref{algorithm:p-levels structured 2}, respectively.
  \State Compute $z = z_H + z_T$ which is the product of $A$ and $v$.
  \end{algorithmic}
\end{algorithm}

\begin{algorithm}
  \caption{$p$ levels sparse matrix-vector product}
    \label{algorithm:p-levels structured 6}
  \begin{algorithmic}[1]
\State Express $A$ by its entries $A = (a_{ij})$.
\State Express $v$ by a column vector 
\begin{multline*}
v=(v_{1,1},\dots, v_{1,N_1}, v_{2,1},\dots, v_{2,N_1},\dots, \\
v_{n_1,1},\dots, v_{n_1,N_1})^{\mathsf{T}}.
\end{multline*}
\State Express $\varphi$ as $(\varphi_1,\dots, \varphi_r)^{\mathsf{T}}$ where $\varphi_j$ is a linear functional on $X_2\circledast \cdots \circledast X_p$ and $r = \prod_{i=2}^p \dim(X_i)$. 
\State Express $\psi$ as $(\psi_1,\dots, \psi_{r})$ where $\psi_j$ is a linear functional on $\mathbb{C}^{N_1}$.
\State Denote $v_i = (v_{i1},\dots, v_{i,N_1})^\mathsf{T}, i=1,\dots,n_1$.
\For{$s=1,\dots,n_1$}
\For{$t=1,\dots,r$}
\State Compute
\[
w_{s,t} = \varphi_t(B) \sum_{(k,s)\not\in \Omega} a_{ks}\psi_t(v_k)
\]
\EndFor
\EndFor
\State Compute 
\[
(z_{ij}) = (I \circledast \vartheta)(w_{st}).
\]
\State $(z_{1,1},\dots, z_{1,n_1},z_{2,1},\dots, z_{2,n_1},\dots, z_{N_1,1},\dots,z_{N_1,n_1})^\mathsf{T}$ is $(A\circledast B)v$.
  \end{algorithmic}
\end{algorithm}

\section{Bilinear complexity}\label{sec:bc}

As we have shown in \cite{YL}, all $11$ algorithms presented in this article are of optimal bilinear complexity, i.e., requires a minimum number of multiplications. We give the multiplication counts below.
\begin{enumerate}[\upshape (i)]
\item  Algorithm~\ref{algorithm:circulant} for $n\times n$ circulant matrix-vector product costs $n$ multiplications (from the computation of $\tilde{z}$; note that the other multiplications in the algorithm are scalar multiplications and do not count towards bilinear complexity).

\item Algorithms~\ref{algorithm:Toeplitz} and \ref{algorithm:Hankel} for $n\times n$ Toeplitz/Hankel matrix-vector products each costs $2n-1$ multiplications (from the computation of $\tilde{z}$; by our special choice of $b$ we saved one multiplication).
 
\item Algorithm~\ref{algorithm:Symmetric}  for $n\times n$ symmetric matrix-vector product costs $\binom{n+1}{2}$ multiplications (each $w_k$ costs $2[n-2(k-1)]-1$ multiplications and so the total number of multiplications is $\binom{n+1}{2}$).

\item An $N\times N$ $p$ levels structured matrix-vector product costs $\prod_{i=1}^p \dim X_i$ multiplications. Let $r = \prod_{i=2}^p \dim(X_i)$.

\item  Algorithm~\ref{algorithm:p-levels structured 1} costs $n_1 r$ multiplications (each $w_{st}$ costs one multiplication; note that computation of the coefficient $\tilde{a}_s\varphi_t(B)$ does not cost any multiplication as $\tilde{a}_s\varphi_t(B)$ is a linear combination of the entries of $A\circledast B$).
\item Algorithms~\ref{algorithm:p-levels structured 2} and Algorithm~\ref{algorithm:p-levels structured 3} each costs $(2n_1-1)r$ multiplications.
\item Algorithm~\ref{algorithm:p-levels structured 4} costs $\binom{n+1}{2}r$ multiplications.
\item Algorithm~\ref{algorithm:p-levels structured 5} costs $(4n-3)r$ multiplications.
\item Algorithm~\ref{algorithm:p-levels structured 6} costs $\#\Omega \times r$ multiplications.
\end{enumerate}

\section*{Acknowledgment}

The authors thank  Nikos Pitsianis and Xiaobai Sun for helpful discussions.
LHL and KY are supported by AFOSR FA9550-13-1-0133, DARPA D15AP00109, NSF IIS 1546413, DMS 1209136, and DMS 1057064. In addition, KY's work is also partially supported by NSF CCF 1017760.

\end{document}